\documentclass[11pt,a4paper,twoside]{article}
\usepackage[UKenglish]{babel}
\usepackage[T1]{fontenc}
\usepackage[utf8x]{inputenc}
\usepackage{latexsym,amsfonts,amsmath,amsthm,amssymb,mathrsfs}
\usepackage{fullpage}
\usepackage{xcolor,bm, bbm}
\usepackage{paralist}
\usepackage{todonotes}
\usepackage{dsfont}
\usepackage{float}
\usepackage{caption}
\usepackage{xfrac}

\usepackage{graphicx}

\DeclareMathOperator{\R}{\mathbb{R}}
\DeclareMathOperator{\N}{\mathbb{N}}

\DeclareMathOperator{\T}{\mathbb{T}}

\DeclareMathOperator{\Var}{Var}
\DeclareMathOperator{\lcm}{lcm}

\newcommand\bq{{\bm{q}}}
\newcommand\br{{\bm{r}}}
\newcommand\Aq{A_{\bq}}
\newcommand\Ar{A_{\br}}

\newtheorem{thm}{Theorem}

\newtheorem{lem}{Lemma}

\newtheorem{proposition}[lem]{Proposition}

\newtheorem*{rem*}{Remark}
\title{\bf Multiplicative D-S}
\medskip

\title{Quantitative inhomogeneous Diophantine approximation for systems of linear forms}
\author{Manuel Hauke}
\date{\today}

\begin{document}

\maketitle
\begin{abstract}
    The inhomogeneous Khintchine–Groshev Theorem is a classical generalization of Khintchine’s Theorem in Diophantine approximation, by approximating points in $\mathbb{R}^m$ by systems of linear forms in $n$ variables. Analogous to the question considered by Duffin and Schaeffer for Khintchine's Theorem (which is the case $m = n = 1$), the question arises for which $m,n$ the monotonicity can be safely removed. If $m = n = 1$, it is known that monotonicity is needed. Recently, Allen and Ramírez showed that for $mn \geq 3$, the monotonicity assumption is unnecessary, conjecturing this to also hold when $mn = 2$. In this article, we confirm this conjecture for the case $(m,n)=(1,2)$ whenever the inhomogeneous parameter is a non-Liouville irrational number. Furthermore, under mild assumptions on the approximation function, we show an asymptotic formula (with almost square-root cancellation), which is not even known for homogeneous approximation. The proof makes use of refined overlap estimates in the 1-dimensional setting, which may have other applications including the inhomogeneous Duffin-Schaeffer conjecture.
\end{abstract}

\section{Introduction and statement of results}

\subsection*{Khintchine-type results}

Khintchine's Theorem is one of the cornerstones in metric Diophantine approximation. It tells us that
for any function $\psi: \N \to [0,\infty)$ that is monotonically decreasing, the Lebesgue measure of the set

\[\{\alpha \in [0,1): \lVert q\alpha \rVert \leq \psi(q) \text{ for infinitely many }q \in \N\}\]
is equal to $0$ resp. $1$, whenever $\sum_{q \in \N} \psi(q)$ converges resp. diverges. 
Here, $\lVert \cdot \rVert$ denotes the distance to the next integer.
Khintchine's Theorem can be generalized to approximating points in arbitrary dimensions by systems of linear forms and allowing an inhomogeneous parameter $\bm{\gamma}$. With these generalizations, this is known as the inhomogeneous Khintchine-Groshev Theorem
\cite{groshev}.\\

\noindent {\bf Inhomogeneous Khintchine-Groshev:} {\it Let 
$m,n \in \mathbb{N}, \bm{\gamma} \in \mathbb{Z}^m$,  and $\psi: \N \to [0,\infty)$ \textbf{monotonically decreasing}. 
Writing $\lvert \bq \rvert:= \max_{1 \leq i \leq n} \lvert q_i \rvert$ and defining

\[ A_{n,m}^{\psi,\gamma} :=
\left\{\bm{\alpha} \in \T^{n \times m}: \lvert \bq\bm{\alpha} + \bm{p} -\bm{\gamma} \rvert \leq \psi(\lvert \bm{q}\rvert) \text{ for infinitely many } (\bm{q},\bm{p}) \in \mathbb{Z}^n \times \mathbb{Z}^m \right\},
\]
we have that

\[
\lambda_{nm}\left(A_{n,m}^{\psi,\gamma}\right) = \begin{cases}
0 &\text { if }  \sum\limits_{q \in \N} q^{n-1}\psi(q)^m < \infty,\\
    1 &\text{ if } \sum\limits_{q \in \N} q^{n-1}\psi(q)^m = \infty,
\end{cases}
\]
where $\lambda_{nm}$ denotes the $nm$-dimensional Lebesgue measure.\\
}
\\
Duffin and Schaeffer \cite{ds} showed that in the case $m = n = 1$ and $\gamma = 0$, the monotonicity assumption in the Khintchine-Groshev Theorem is indeed necessary, a result generalized by Ramírez \cite{ramirez_counter} to arbitrary $\gamma$. After contributions in higher dimensions by Schmidt \cite{schmidt} and Gallagher \cite{gallagherII}, it was shown by Beresnevich and Velani \cite{bv} that for $\bm{\gamma} = 0$, the monotonicity assumption can be safely removed whenever $mn \geq 2$, a result that completed the characterization in the homogeneous setting. For $\bm{\gamma} \neq 0$, the first result is due to Schmidt \cite{schmidt} \footnote{Although only stated in Schmidt's paper for the homogeneous setup, the proof can be adapted to hold in a very general setup that includes the inhomogeneous case. See Sprind\v{z}uk \cite{sprindzuk} for details.} who showed that the monotonicity can be removed when $n \geq 3, m \geq 1$. There was not much progress in this question until recently: Yu \cite{yu} proved the same for $m \geq 3, n\geq 1$ and Allen and Ramírez \cite{ar} showed that $mn \geq 3$ suffices. They conjectured this to be true also for the case $mn = 2$, which would be in accordance to the homogeneous setting:\\

\noindent {\bf Conjecture (Allen, Ramírez).}
{\it 
The monotonicity condition in the inhomogeneous Khintchine-Groshev Theorem can be removed also in cases $(m,n) = (2,1)$ and $(m,n) = (1,2)$.
}\\

Assuming some extra divergence conditions in the vein of 
\cite{az,bv_extra,hpv}, Yu \cite{yu} respectively Allen and Ramírez \cite{ar} showed the conjecture to be true. However, the conjecture is significantly harder to solve without assuming extra divergence and as such, the conjecture is completely open.\\

In this paper we make progress towards proving the conjecture: As a first result in this article, we show that if $(m,n) = (1,2)$ and $\gamma$ is a non-Liouville irrational number, the monotonicity condition in the inhomogeneous Khintchine-Groshev Theorem can be removed, without imposing any condition on $\psi$.

\begin{thm}\label{khintchine_groshev}
Let $\psi: \N \to [0,\infty)$ be an arbitrary function and 
$\gamma$ a non-Lioville irrational. Then we have

\[
\lambda_{2}\left(A_{2,1}^{\psi,\gamma}\right) = \begin{cases}
0 &\text { if } \sum\limits_{q \in \N} q\psi(q) < \infty,\\
    1 &\text{ if } \sum\limits_{q \in \N} q\psi(q) = \infty.
\end{cases}
\]
\end{thm}

We note that from a metrical point of view, the condition of being a non-Liouville irrational is a rather weak restriction: Although Liouville numbers form an uncountable set, they have Hausdorff dimension $0$ (and thus in particular, Lebesgue measure $0$). For a proof of these results we refer the reader to standard literature such as \cite{oxtoby}.\\

Using the framework of Allen and Beresnevich \cite{ab_hausdorff}, which generalizes the mass transference principle established by Beresnevich and Velani \cite{bv_mass}, we obtain as immediate corollary the following Hausdorff analogue to Theorem \ref{khintchine_groshev}: 
Defining

\[
t := \inf\left\{ s \in [1,\infty): \sum_{q \in \N} q^2 \left(\frac{\psi(q)}{q}\right)^{s-1}  < \infty\right\},
\]
we obtain that for non-Liouville irrational $\gamma$, the Hausdorff dimension of ${A_{2,1}^{\psi,\gamma}}$ satisfies \[\dim_{\mathcal{H}}({A_{2,1}^{\psi,\gamma}}) = \min\{t,2\}.\]

For the proof of this fact and more detailed definitions, we refer the reader to \cite{ar} and the references therein.

\subsection*{The quantitative theory}

For fixed $\bm{\alpha} \in \T^n \times \T^m, \psi: \N \to [0,1/2), \bm{\gamma} \in \mathbb{Z}^n$ and a height $Q \in \N$, let

\[
\mathcal{N}(\bm{\alpha},Q,\bm{\gamma}) :=
\#\left\{
(\bm{p},\bm{q}) \in \mathbb{Z}^m \times \mathbb{Z}^n:
\lvert \bm{q} \rvert \leq Q: \left\lvert \bm{q}\bm{\alpha} - \bm{p} - \bm{\gamma} \right\rvert \leq \psi(\lvert \bm{q} \rvert)
\right\}.
\]

Assuming to be in the divergent case and an analogue of Theorem \ref{khintchine_groshev} to be known, then almost surely, $\lim_{Q \to \infty} \mathcal{N}(\bm{\alpha},Q,\bm{\gamma}) = \infty$. Thus, an obvious question to ask is about the existence of an asymptotic formula. If $\psi$ is monotonic and $\bm{\gamma} = 0$, Gallagher \cite{gallagher} proved that in the case $n = 1, m\geq 1$, for almost every $\bm{\alpha} \in \T^m$ and any $\delta > 0$,

\begin{equation}\label{gallagher_asymp}
\mathcal{N}(\bm{\alpha},Q,\bm{\gamma}) = \Psi(Q) + O_{\delta}\left(\Psi(Q)^{1/2}\log(\Psi(Q)^{3/2 + \delta})\right), \quad Q\to \infty,
\end{equation}
where $\Psi(Q)$ equals the corresponding sum of measures, which is (up to constants depending on $n,m$) of the form
$ \sum\limits_{q  \leq Q} {q^{n-1}}\psi(q)^m$.

Without assuming monotonicity, Schmidt \cite{schmidt} showed that for $\bm{\gamma} = 0$,

\begin{equation}\label{schmidt_asymptotic}
\mathcal{N}(\bm{\alpha},Q,\bm{\gamma}) = \Psi(Q) + O_{\delta}\left(\chi(Q)^{1/2}\log(\chi(Q)^{3/2 + \delta})\right), \quad Q\to \infty,
\end{equation}
where 
\[
\chi(Q) := \sum_{\bq \in \mathbb{Z}^n \setminus \{\bm{0}\}} (\psi(\lvert\bq\rvert))^m \tau(\gcd(\bq)),
\]
with $\tau$ denoting the number of divisors function. Employing some modifications, \eqref{schmidt_asymptotic} can be adapted to also hold in the inhomogeneous setting, see e.g. \cite[Theorems 18 and 19]{sprindzuk}. If $n \geq 3$, Schmidt showed that $\chi(Q) \ll \Psi(Q)$ and thus, \eqref{schmidt_asymptotic} is a meaningful result for $n \geq 3$.
However, as beautifully illustrated in \cite{bv}, for $n = 2$ and any $m \in \N$, it is possible to construct $\psi$ such that the ``error term'' dominates the ``main term'', thus the statement \eqref{schmidt_asymptotic} becomes trivial if $\psi$ is chosen in such a way.\\

Here, we show that under mild conditions on $\psi$ and the same non-Liouville condition as in Theorem \ref{khintchine_groshev}, we recover \eqref{gallagher_asymp} for $(m,n) = (1,2)$ without assuming monotonicity and thus improving upon Schmidt's result by getting rid of the divisor function in \eqref{schmidt_asymptotic}:

\begin{thm}\label{main_thm_new}
    Let $\gamma$ be an irrational non-Liouville number, $\varepsilon > 0$ and $\psi: \N \to [0,1/2]$ such that
    \[
    \sum_{q \in \N} q\psi(q) = \infty, \quad \psi(q) = O(1/q^{\varepsilon}).
    \]
    Then for any $\delta > 0$ and almost every $\bm{\alpha} \in \T^2$, we have that
\[
\mathcal{N}(\bm{\alpha},Q,\bm{\gamma}) = 
\Psi(Q) + O\left(\Psi(Q)^{1/2}\log(\Psi(Q)^{3/2 + \delta})\right),\quad Q \to \infty,
\]
where $\Psi(Q) := 16\sum\limits_{q \leq Q}q\psi(q) + 8\sum\limits_{q \leq Q}\psi(q)$.
\end{thm}

\subsection*{Remarks on the theorems and further research}
\begin{itemize}
\item
We note that Theorem \ref{main_thm_new} is somehow surprising: Usually, the homogeneous case (that is, $\gamma = 0$) is much better understood and easier to control. However, to the best of the author's knowledge, there is no asymptotic result for the homogeneous setting. The best result is due to Beresnevich and Velani \cite{bv} who showed an analogue of Theorem \ref{khintchine_groshev} for $(n,m) = (2,1)$, but no asymptotic. The method of their proof, which is standard in metric homogeneous problems, removes non-coprime contributions, that is, only considering those $(p,q_1,q_2) \in \mathbb{Z}^3$ with
$\left\lvert q_1\alpha_1 + q_2\alpha_2 - p \right\rvert \leq \psi(\lvert(q_1,q_2)\rvert)$
such that $\gcd(q_1,q_2,p) = 1$. By doing this, one loses a constant factor, which is enough to establish the Khintchine-type result from Theorem \ref{khintchine_groshev}.
It might be possible to work out the asymptotic for the number of \textit{coprime} good approximations similarly to \cite{abh}, but clearly this approach does not provide any hope of showing an asymptotic formula for \textit{all} $(q_1,q_2,p)$ as it is obtained in Theorem \ref{main_thm_new}. However, without imposing an extra condition or removing some elements, the variance of the set system (in the homogeneous case) can be much larger (which is also reflected in the counterexamples \cite{ds} and \cite{ramirez_counter} when $m = n = 1$) and thus, the standard techniques are impossible to use in order to prove an analogue of Theorem \ref{main_thm_new} in the homogeneous setup.
Thus Theorem \ref{main_thm_new} is one of the rare occasions where results in the inhomogeneous case are known that are still out of reach for the homogeneous counterpart.\\
For Liouville numbers (and even for rational $\gamma \neq 0$), the statement is still wide-open: The error terms we accumulate here become too large and we could only avoid them by assuming an extra divergence factor similar to \cite{ar,yu}. Thus, probably some form of 
coprimality condition needs to be imposed in these cases. However, similarly to the even more challenging $1$-dimensional case, it remains unclear which condition should be imposed to establish Theorem \ref{khintchine_groshev}, let alone a quantification in the form of Theorem \ref{main_thm_new}.
\item  If $n = 2, m \geq 2$, an analogue of Theorem \ref{main_thm_new} holds without any assumption on the inhomogeneous parameter $\gamma$ (in fact, the inhomogeneous shift could be different for each $\bm{q}$ as in the general setup considered in \cite{ar}.) If $n \geq 3$ and $m\geq 1$, even the assumption $\psi(q) = O(1/q^{\varepsilon})$ can be dropped, which then coincides exactly with the result that can be recovered from Schmidt's asymptotic formula \eqref{schmidt_asymptotic}. We provide a sketch of the necessary adaptation of the proof of Theorem \ref{main_thm_new} to the cases $(n,m) = (2,2)$ and $n \geq 3$
at the end of this article. Note that this recovers the (already known) Khintchine-type result for $n \geq 2, m\geq 2$ since replacing $\psi$ with  $\tilde{\psi}(q) := \min\{ \frac{1}{q^{\varepsilon}},\psi(q)\}$ does not alter the divergence property of
$\sum_{q \in \N} q^{n-1}\psi(q)^m.$
\item Schmidt \cite{schmidt} respectively Sprind\v{z}uk \cite{sprindzuk} actually stated \eqref{schmidt_asymptotic} in a more flexible way by allowing $\psi$ to be a multivariate function (that is, $\psi$ is not only depending on $\lvert\bm{q}\rvert$, but on $\bm{q}$ itself). In this setup, monotonicity is definitely needed for any $n$ and any $\gamma$: Considering functions $\psi$ that are only supported on vectors where all except the first coordinate vanish, the counterexamples of Duffin and Schaeffer \cite{ds} and Ramírez \cite{ramirez_counter} show the necessity of monotonicity in this setting. Thus, there is no chance to generalize Theorem \ref{main_thm_new} in that way.

    \item
The error term in Theorem \ref{main_thm_new} is, up to logarithmic factors, best possible when considering a mean-variance argument - the variance will turn out to be bounded (up to constants) by its diagonal contribution, which is an unavoidable error term. The method of proof shows that the implied constant in Theorem \ref{main_thm_new} is indeed absolute, not depending on $\bm{\alpha},\gamma,\delta,\psi$, as long as $Q$ is large enough in  terms of these parameters. To be more precise, for any admissible $\bm{\alpha},\gamma,\delta,\psi$ there exists an ineffective constant $K(\bm{\alpha},\gamma,\delta,\psi)$ and an absolute effective constant $C >0$ such that 

\[\lvert \mathcal{N}(\bm{\alpha},Q) - 
\Psi(Q) \rvert \leq K(\bm{\alpha},\gamma,\delta,\psi) + C\Psi(Q)^{1/2}\log(\Psi(Q)^{3/2 + \delta}).\]

\end{itemize}

We conclude this introductory section by making some heuristic remarks on the kind of problems considered in this article: The higher the dimensions of $n$ and $m$ in the Khintchine-Groshev-Theorem, the easier it becomes to get an almost sure behaviour or even an asymptotic formula - an observation that was formalized in an elegant way by Allen and Ramírez \cite{ar}. One reason for this can be seen from the fact that the random events of $\bq \in \mathbb{Z}^n$ respectively $\br \in \mathbb{Z}^n$ contributing something to $\mathcal{N}(\bm{\alpha},Q)$ for non-parallel vectors $\bq \nparallel \br$ are independent (see Proposition \ref{nonparallel}). The higher the dimension $n$, the smaller the proportion of parallel vectors where non-immediate bounds on the degree of independence have to be found, which makes the task of establishing useful variance estimates much easier. \\
In the case of parallel vectors, the problem is equivalent to treating $m$-dimensional overlap estimates, which, if $m$ is large does not cause too many problems. However if $m = 2$ and even more so if $m = 1$, the overlap estimates (especially in the inhomogeneous setting) are far from completely understood.
The main difficulty is, even more so in the inhomogeneous setting, to remove (almost) identical overlaps. To establish Khintchine-type results in the homogeneous world, this can be overcome  as already mentioned, by restricting to coprime solutions. This removes the possibility of identical overlaps completely, which is, amidst ingenious further ingredients due to Koukouloupolos and Maynard \cite{km}, the reason that the Duffin-Schaeffer conjecture is true. In the inhomogeneous irrational setting, an exact overlap cannot occur, a fact that arises from
$q\gamma \notin \mathbb{Z}$ for any $q \in \mathbb{Z} \setminus \{0\}$, and the property of $\gamma$ being a non-Liouville irrational gives in the two-dimensional setting considered here enough information about $\lVert q\gamma \rVert$ staying away from zero to establish a refined $1$-dimensional overlap estimate (Lemma \ref{overlap_estim_new}), which provides the new key ingredient to establish the main results in this article. From a technical point of view, 
the approach of Allen and Ramírez \cite{ar} makes use of 

\[\frac{1}{q^{n+m-1}}\sum\limits_{\substack{r\leq q}}
\gcd(q,r)^{m+n-1} \ll 1,\] which is true whenever $m+n \geq 4$. By the assumption of non-Liouville and assuming $\psi(q) = O(q^{-\varepsilon})$, our improved overlap estimate leads to
\[\frac{1}{q^{n+m-1}}\sum\limits_{\substack{r\leq q\\\gcd(q,r) < q^{1-\delta_{\varepsilon}}}}
\gcd(q,r)^{m+n-1} \ll 1,\]
which is true when $m+n-1 \geq 3$. 
Note that $\psi(q) = O(q^{-\varepsilon})$ is no restriction when considering $(n,m) = (2,1)$ since $\sum\limits_{q \in \N} q\min\left\{\psi(q),\frac{1}{q^{\varepsilon}}\right\}$ still diverges if the original sum $\sum\limits_{q \in \N} q\psi(q)$ does. This is not true in the $(n,m) = (1,2)$ case since 
it can happen that $\sum\limits_{q \in \N} \psi(q)^2$ diverges, but $\sum\limits_{q \in \N} \left(\min\left\{\psi(q),\frac{1}{q^{\varepsilon}}\right\}\right)^2$ converges.
    
\section{Notation and Preliminary results}

\subsection*{Notation}
For a set $A$, we write the indicator function as $\mathds{1}_A(\cdot)$.
We use the standard $O$- and $o$-symbols as well as the Vinagradov notation
$f \ll g :\Leftrightarrow f = O(g)$.
We write $\T^k = \left(\sfrac{\R}{\mathbb{Z}}\right)^k$ and denote by $\lambda_k$ the $k$-dimensional normalized Haar measure on $\T^k$. Given a vector $\bq = (q_1,\ldots q_n) \in \mathbb{Z}^n$, we write $\lvert \bq \rvert := \max\limits_{1 \leq i \leq n}\lvert q_i\rvert$ and 
$\gcd(\bq) = \gcd(q_1,\ldots,q_n)$. To emphasize the difference between vectors and scalars, we write vectors as $\bq$ and scalars as $q$.

For $\bq,\br \in \mathbb{Z}^2$, we write $\bq \parallel \br$ whenever there
exist $n,m \in \mathbb{Z}$ such that $n\bq = m\br$. If this is not the case, we write
$\bq \nparallel \br$.
\\

A number $\gamma$ is called \textit{Liouville number} if for any $\eta \in \N$, there exists $q \in \N$ with $q \geq 2$ such that 
$
\left\lVert  q\gamma \right\rVert  \leq \frac{1}{q^{\eta}}$.
Thus, for every non-Liouville number $\gamma$ there exists $\eta = \eta(\gamma) \in \mathbb{N}$ and $c = c(\gamma) > 0$ such that for all $q \neq 0$,

\begin{equation}\label{non_liouville}
\left\lVert  q\gamma \right\rVert \geq \frac{1}{cq^{\eta}}.
\end{equation}
For notational conveniencs, we extend any $\psi: \mathbb{N} \to \R$ to $\mathbb{Z}^n$ by defining $\psi(\bq):= \psi(\lvert \bq \rvert )$.

\subsection*{Preliminaries}
For fixed $\psi: \mathbb{N} \to \R$, $\gamma \in \R$, we define the sets 
\[A_{\bq} := A_{\bq,\psi,\gamma}
:= \left\{\bm{\alpha} \in \T^2: \lVert \bq\bm{\alpha} - \gamma \rVert \leq \psi( \bq )
\right\} \subseteq \T^2.
\]

The statements of Theorems  \ref{khintchine_groshev} and \ref{main_thm_new} can be rephrased in terms of $\Aq$ and will be proved by a mean-variance argument, which is the usual approach in metric Diophantine approximation. We will apply it in the formulation of \cite[Lemma 1.5]{harm}.

  \begin{lem}\label{harman}
        Let $X$ be a measure space with measure $\mu$ such that $0 < \mu(X) < \infty$.
        Let $(f_k(x))_{k \in \N}$ be a sequence of non-negative $\mu$-measurable functions and let $f_k, \varphi_k$ be real numbers such that for all $k \in \N$,
        \[
        0 \leq f_k \leq \varphi_k.
        \]
        Write $\Phi(N) = \sum_{k = 1}^{N} \varphi_k$ and suppose that $\Phi(N) \to \infty$ as $N \to \infty$. Suppose that for arbitrary integers $1 \leq m < n$ we have
        \begin{equation}\label{local_variance}
        \int_{X} \left(\sum_{m \leq k < n} (f_k(x) - f_k)\right)^2 \mathrm{d}\mu \leq K\sum_{m \leq k < n}\varphi_k
            \end{equation}
        for an absolute constant $K$. Then for any given $\delta > 0$ and almost all $x \in X$, we have as $N \to \infty$,

        \[\sum_{k = 1}^N f_k(x) = 
        \sum_{k=1}^N f_k + O\left(\Phi^{1/2}(N) \left(\log(\Phi(N) +2)\right)^{3/2 + \delta} + \max_{1 \leq k \leq N}f_k\right).
        \]
    \end{lem}

In our setting, we will have $f_k(x) = \mathds{1}_{A_{\bq_k}}(x)$
for some $\bq_k \in \mathbb{Z}^2$ and $f_k$ will denote the corresponding expected value. 
Thus in view of Lemma \ref{harman}, we need to compute the expected value of the sum, that is, the sum of the measures.

\begin{proposition}\label{sum_of_measures}
For any function $\psi: \N \to [0,1/2]$ and any $\gamma \in \R$, we have

\[
\sum_{q \leq Q}\sum_{\substack{\bq \in \mathbb{Z}^2\\ \lvert \bq\rvert = q}}
\lambda_2(\Aq) =
16\sum_{q \leq Q}q\psi(q) + 8\sum_{q \leq Q}\psi(q).
\]
\end{proposition}

\begin{proof}
An elementary computation shows that
\[
\lambda_2(\Aq) = 2\psi(\bq), \quad
\sum_{q \leq Q}\sum_{\substack{\bq, \in \mathbb{Z}^2\\ \lvert \bq\rvert = q}}1
 = 8q + 4,
\]
which proves the statement.
\end{proof}

Inequality \eqref{local_variance}, which can be interpreted as a variance estimate, is from a technical perspective, the key ingredient to proving theorems in metric Diophantine approximation. We state 
the version needed for Theorems \ref{khintchine_groshev} and \ref{main_thm_new} in the following lemma.

 \begin{lem}\label{sneaky_order}
        Let $\prec$ denote a total order of $\mathbb{Z}^2$ that satisfies
        \[
        \lvert \br \rvert < \lvert \bq \rvert \Rightarrow \br \prec \bq.
        \]
        Let $\gamma$ be a non-Liouville irrational and $\psi: \N \to [0,1/2]$ with $\psi(q) = O(q^{-\varepsilon})$.
        Then for any $\bm{u} \prec \bm{v}$ 
        and $\psi,\gamma$ as before, we have that

        \begin{equation}\label{fine_scale_variance}
\sum_{\bm{u} \preceq \bq,\br \preceq \bm{v}}
\lambda_2\left(\Aq \cap \Ar\right) 
- \left(\sum_{\bm{u} \preceq \bq,\br \preceq \bm{v}}\lambda_2\left(\Aq\right)\right)^2 \ll
\sum_{\bm{u} \preceq \bq,\br \preceq \bm{v}}\lambda_2\left(\Aq\right),
        \end{equation}
        with the implied constant being independent of $\bm{u},\bm{v}$.
    \end{lem}

The proof of Lemma \ref{sneaky_order} contains the main difficulty in obtaining the results of this paper, thus we postpone it to Section \ref{lemma_proof}. Now we show how Theorems \ref{khintchine_groshev} and \ref{main_thm_new} follow from it.

    \begin{proof}[Proof of Theorem \ref{main_thm_new}, assuming Lemma \ref{sneaky_order}]
    We apply Lemma \ref{harman} in the following way: Let $X = \T^2, \mu = \lambda_2$.
    Using any total order as in Lemma \ref{sneaky_order}, we order all of $\mathbb{Z}^2$ by
    \[
\bq_1 \prec \bq_2 \prec \bq_3 \ldots
    \]
    and define
    $f_k(x) = \mathds{1}_{A_{\bq_k}}(x), f_k = \varphi_k = \lambda_2\left(A_{\bq_k}\right) \asymp 
    \psi(\bq_k)$.  Applying Lemma \ref{sneaky_order}, we see that \eqref{local_variance} is satisfied. By Proposition \ref{sum_of_measures} and since $f_k \leq 1$, the statement follows from Lemma \ref{harman}.
    \end{proof}

    \begin{rem*}
        Actually, we are proving something slightly stronger here: For any total order as in Lemma \ref{sneaky_order} and any $\bm{u} \in \mathbb{Z}^2$, we obtain an asymptotic formula (with the corresponding error term) for $\{\bq \in \mathbb{Z}^2: \bq \prec \bm{u}, \lVert \bq\alpha - \gamma \rVert \leq \psi(\lvert \bq\rvert)\}$ as $\bm{u} \to \infty$ (in the sense of the total order). Theorem 2 is the special case where $\bm{u}$ is of the form
        $\bm{u} \prec \bm{q} \Rightarrow \lvert \bm{u}\rvert < \lvert\bm{q}\rvert$. However, to prove this special case, we cannot only work with such $\bm{u}$ since grouping all $A_{\bm{q}}$ with same norm $q$, we have (in the notation of Lemma \ref{harman}) $f_q = q\psi(q)$, which is in general unbounded.
    \end{rem*}

    \begin{proof}[Proof of Theorem \ref{khintchine_groshev} assuming Theorem \ref{main_thm_new}]
        If $\sum_{q \in \N} q\psi(q)< \infty$, the statement follows by applying the convergence Borel-Cantelli Lemma. Thus we can assume $\sum_{q \in \N} q\psi(q)= \infty$.
        Furthermore, we can assume $\psi(q) \leq \frac{1}{2q}$ for any $q \in \N$ since otherwise, we replace $\psi$ by $\tilde{\psi}(q) := \min\{\frac{1}{2q},\psi(q)\}$, which still satisfies $\sum_{q \in \N} q\tilde{\psi}(q)= \infty$. Therefore, we can apply Theorem \ref{main_thm_new} and the statement follows.
    \end{proof}

\section{Proof of Lemma \ref{sneaky_order}}\label{lemma_proof}
\subsection*{Overlap estimates for parallel vectors}

In order to prove Lemma \ref{sneaky_order}, we have to control the degree of dependence of sets
$\Aq,\Ar$. A main advantage of linear forms with at least two variables in comparison to the 1-dimensional setting is that for non-parallel vectors (which is the ``generic'' case for dimension $n \geq 2$), the sets
$\Aq, \Ar$ are independent:

\begin{proposition}\label{nonparallel}
    Let $\bq,\br \in \mathbb{Z}^2$ such that $\bq \nparallel \br$ and $\gamma \in \R$ arbitrary. Then 
      \[
    \lambda_2(\Aq \cap \Ar)
    =  \lambda_2\left(\Aq\right) \lambda_2\left(\Ar\right).
    \]
\end{proposition}

\begin{proof}
 This follows from a well-known more general statement, for details see e.g. \cite[Lemma 9, Chapter 1]{sprindzuk}.
\end{proof}

To treat the case of parallel vectors, we will need a more delicate analysis. 
Since we assumed $\psi(q) = O(q^{-\varepsilon})$ and $\gamma$ non-Liouville, there exists $\eta = \eta(\gamma) \in \N,C,c > 0$ such that for all $q \in \mathbb{Z}\setminus\{0\}$, we have
\begin{equation}\label{assumptions}
\psi(q) \leq \min\left\{\frac{C}{q^{\varepsilon}},1/2\right\}, \quad
\lVert q\gamma\rVert \geq \frac{1}{cq^{\eta}}.
\end{equation}
Thus the following lemma is key to establish the variance estimate:

\begin{lem}\label{overlap_estim_new}
    Let $r,q \in \N$ with $r < q$ and assume \eqref{assumptions} holds. Further, let $\bq\parallel \br$ be vectors in $\mathbb{Z}^2$ with $\lvert \bq \rvert = q, \lvert \br \rvert = r, d = \gcd(\bq), e = \gcd(\br)$. 
Then
\[
\lambda_2(\Aq \cap \Ar)
\leq \begin{cases}
    0 &\text{ if } r > d^{\frac{\eta+1}{\varepsilon}}(2Cc)^{\frac{1}{\varepsilon}},\\
    \lambda_2(\Aq) \lambda_2(\Ar) + 4\frac{\psi(q)}{d}\gcd(d,e) &\text{ otherwise }.
\end{cases}
\]
\end{lem}

\begin{proof}
    Since $\bq \parallel \br$ and $d = \gcd(\bq), e = \gcd(\br)$, we have
$
\bq = \pm d\bm{p}, \br = \pm e \bm{p}$,
where $\bm{p}$ is a primitive vector (that is, $\bm{p} = (p_1,p_2)$ with $\gcd(p_1,p_2) = 1$).
We assume without loss of generality that $\bq = d\bm{p}, \br = e\bm{p}$.
Observe that (compare e.g. \cite[Lemma 5]{ar}) 
\[
\Aq \cap \Ar
= T_{\bm{p}}^{-1}\left( A(d,\psi(q)) \cap A(e,\psi(r))\right)
\]
where $T_{\bm{p}}: \T^2 \to \T: \bm{x} \mapsto \bm{p}\bm{x}$
and
$
A(d,t) := \{\alpha \in \T: \lVert d\alpha - \gamma\rVert \leq t\}
$
denotes $1$-dimensional sets.
Since $T_{\bm{p}}$ is measure-preserving, we thus have 

\[
\lambda_2(\Aq) =  \lambda_1(A(d,\psi(q)), \quad \lambda_2(\Ar) =  \lambda_1(A(e,\psi(r))
\]
as well as
\[
\lambda_2(\Aq \cap \Ar)= \lambda_1(A(d,\psi(q)) \cap A(e,\psi(r))).
\]
Thus we are left to obtain an overlap estimate in dimension $1$. Denoting

\[\begin{split}
\delta &:= \min\left\{\frac{\psi(q)}{d}, \frac{\psi(r)}{e}\right\},\quad \Delta := \max\left\{\frac{\psi(q)}{d}, \frac{\psi(r)}{e}\right\},
\end{split}
\]
and defining the piecewise linear function 

\begin{equation}\label{weight_fct}
w(y) := \begin{cases}
    2\delta &\text{ if } 0 \leq \lvert y \rvert \leq \Delta - \delta,\\
    \Delta + \delta - y &\text{ if } \Delta - \delta < \lvert y \rvert \leq \Delta + \delta,\\
    0 &\text{ otherwise,}
\end{cases}
\end{equation}
we see that (compare \cite[Lemma 5]{abh} or \cite[Proposition 7]{lonely})
\[
\lambda_1(A(d,\psi(q)) \cap A(e,\psi(r))) 
= \sum_{a \in \mathbb{Z}_d}\sum_{b \in \mathbb{Z}_e}
w\left(\frac{a}{d} - \frac{b}{e} - \left(\frac{\gamma}{d} - \frac{\gamma}{e}\right)\right).
\]
Since every $c \in \mathbb{Z}_{\lcm(d,e)}$ has exactly 
$\gcd(d,e)$ many solutions $(a,b) \in \mathbb{Z}_d \times \mathbb{Z}_e$ to the equation
\[
\frac{a}{d} - \frac{b}{e} = \frac{c}{\lcm(q,r)}, 
\]
we obtain
\[
\lambda_1(A(d,\psi(q)) \cap A(e,\psi(r))) 
= \gcd(d,e)\sum_{c \in \mathbb{Z}_{\lcm(d,e)}}
w\left(\frac{c - \gamma\left(\frac{d-e}{\gcd(d,e)}\right)}{\lcm(d,e)}\right).
\]

We apply an argument used by Schmidt \cite{schmidt} (see \cite[Lemma 4.4]{harm} for a detailed study), by removing the $c$ close to the origin and using the monotonicity of $w$. To do this formally, let 
\[c_1 := 
\left\lfloor  \gamma\left(\frac{d-e}{\gcd(d,e)}\right)\right\rfloor, \quad
c_2 := \left\lceil  \gamma\left(\frac{d-e}{\gcd(d,e)}\right)\right\rceil.
\]
Since $w(y) \leq 2\delta$ for all $y \in \R$, we can bound

\[\begin{split}
\gcd(d,e)\sum_{c \in \mathbb{Z}_{\lcm(d,e)}}
w\left(\frac{c - \gamma\left(\frac{d-e}{\gcd(d,e)}\right)}{\lcm(d,e)}\right)
&\leq 4\delta \gcd(d,e) +  \gcd(d,e)\sum_{\substack{c \in \mathbb{Z}\\c \notin \{c_1,c_2\}}}
w\left(\frac{c - \gamma\left(\frac{d-e}{\gcd(d,e)}\right)}{\lcm(d,e)}\right)
\\&\leq 4\delta \gcd(d,e) +  \gcd(d,e)\lcm(d,e)\int_{-\infty}^{\infty} w(x) \,\mathrm{d}x,
\end{split}
\]
where we used that $c_1 - \gamma\left(\frac{d-e}{\gcd(d,e)}\right) < 0 < c_2 - \gamma\left(\frac{d-e}{\gcd(d,e)}\right)$ and $w$ being monotonically increasing on $(-\infty,0)$ and monotonically decreasing on $[0,\infty)$. Since 
\[\gcd(d,e)\lcm(d,e)\int_{-\infty}^{\infty} w(x) \,\mathrm{d}x = 4\gcd(d,e)\lcm(d,e)\delta \Delta =  \lambda_2(\Aq) \lambda_2(\Ar),\] the result follows for arbitrary $r$. If $r$ is large, we can do better under the assumption of \eqref{assumptions}, by showing that for no $c \in \mathbb{Z}$, $\frac{c -  \gamma\left(\frac{d-e}{\gcd(d,e)}\right)}{\lcm(d,e)}$ lies in the support of $w$,
or in other words,
\begin{equation}\label{no_close_integer}
\left\lVert \gamma\left(\frac{d-e}{\gcd(d,e)}\right)\right\rVert > (\delta + \Delta)\lcm(d,e).
\end{equation}
Since by \eqref{assumptions}, we have
$\psi(m) \leq \frac{C}{m^{\varepsilon}}$ and hence
\[
(\delta + \Delta)\lcm(d,e) \leq 
2 de
\max\left\{\frac{C}{r^{\varepsilon}e}, \frac{C}{q^{\varepsilon}d}\right\}
= \frac{2Cd}{r^{\varepsilon}}.
\]
On the other hand, \eqref{assumptions} also implies (note that $d \neq e$ since $q \neq r$),
\[\left\lVert \gamma\left(\frac{d-e}{\gcd(d,e)}\right) \right\rVert 
\geq \frac{1}{cd^{\eta}}.
\]
Thus, for 
$\frac{2Ccd}{r^{\varepsilon}}
< \frac{1}{d^{\eta}}
$ 
we obtain \eqref{no_close_integer}, which concludes the proof.

\end{proof}

\subsection*{Finishing the proof of Lemma \ref{sneaky_order}}

 In this section, we combine the previous results to prove Lemma \ref{sneaky_order}. In order to do so, we first prove the special case where $\bm{u} = (0,0)$ and $\bm{v}$ satisfies
 $\bm{v} \prec \bq \Rightarrow \lvert \bm{v} \rvert < \lvert \bq \rvert$.
 
\begin{lem}\label{variance_estim_new}
  Let $\psi(m) = O(q^{-\varepsilon})$ and $\gamma$ a non-Liouville irrational.
    Then we have

    \[\sum_{q,r \leq Q}
\sum_{\substack{\bq \in \mathbb{Z}^2\\ \lvert \bq\rvert = q}}
\sum_{\substack{\br \in \mathbb{Z}^2\\ \lvert \br\rvert = r}}\lambda_2(\Aq \cap \Ar) 
= \Psi(Q)^2 + O\left(\Psi(Q)\right), \quad Q\to \infty.
\]
\end{lem}
\begin{proof}[Proof of Lemma \ref{variance_estim_new}]
    We first split the sums into parallel vectors and non-parallel ones. By Proposition \ref{nonparallel}, the sets corresponding to non-parallel vectors are perfectly independent, so we are left to consider the parallel cases. We have 

\begin{align} \sum_{q,r \leq Q}
\sum_{\substack{\bq, \in \mathbb{Z}^2\\ \lvert \bq\rvert = q}}
\sum_{\substack{\br \in \mathbb{Z}^2\\ \lvert \br\rvert = r\\ \br \parallel \bq}} \lambda_2(\Aq \cap \Ar)
 = & \sum_{q \leq Q}\sum_{\substack{\bq, \in \mathbb{Z}^2\\ \lvert \bq\rvert = q}}
\sum_{\substack{\br \in \mathbb{Z}^2\\ \lvert \br\rvert = q\\ \br \parallel \bq}} \lambda_2(\Aq \cap \Ar)\label{diag}
\\+
2 & \sum_{q \leq Q}\sum_{\substack{\bq, \in \mathbb{Z}^2\\ \lvert \bq\rvert = q}}
\sum_{r < q}\sum_{\substack{\br \in \mathbb{Z}^2\\ \lvert \br\rvert = r\\ \br \parallel \bq}} \lambda_2(\Aq \cap \Ar).\label{non_diag}
\end{align}

Note that $\br \parallel \bq$ and $\lvert \br \rvert = \lvert \bq \rvert$ implies $\br = \pm \bq$ and thus we can bound the right-hand side of \eqref{diag} by
\begin{equation}\label{diag_contr}
 \sum_{q \leq Q}\sum_{\substack{\bq \in \mathbb{Z}^2\\ \lvert \bq\rvert = q}}
\sum_{\substack{\br \in \mathbb{Z}^2\\ \lvert \br\rvert = q\\ \br \parallel \bq}} \lambda_2(\Aq \cap \Ar)
\leq 2 \sum_{q \leq Q}\sum_{\substack{\bq, \in \mathbb{Z}^2\\ \lvert \bq\rvert = q}}\lambda_2(\Aq) \ll \Psi(Q),
\end{equation}
and we are left to estimate \eqref{non_diag}.
For fixed $q$, we rewrite the contribution of \eqref{non_diag} by

\[
\sum_{\substack{\bq \in \mathbb{Z}^2\\ \lvert \bq\rvert = q}}
\sum_{r < q}\sum_{\substack{\br \in \mathbb{Z}^2\\ \lvert \br\rvert = r\\ \br \parallel \bq}} \lambda_2(\Aq \cap \Ar)
=
\sum_{d\mid q}\sum_{\substack{\bq \in \mathbb{Z}^2\\ \lvert \bq\rvert = q\\
\gcd(\bq) = d}}
\sum_{r < q}
\sum_{\substack{\br \in \mathbb{Z}^2\\ \lvert \br\rvert = r\\ \br \parallel \bq}} \lambda_2(\Aq \cap \Ar).
\]

Let $\bq = d\bm{p}$ be fixed with $\bm{p}$ a primitive vector. If $\br \parallel \bq$ with $\lvert \br \rvert = r$, then $\frac{q}{d} \mid r$ and
$\br = \pm \frac{r}{q/d}\bm{p}$. Thus for fixed $\bq$ such that $\bq = d\bm{p}$, we have

\[\begin{split}&\sum_{r < q}
\sum_{\substack{\br \in \mathbb{Z}^2\\ \lvert \br\rvert = r\\ \br \parallel \bq}} \lambda_2(\Aq \cap \Ar) -
\lambda_2(\Aq)\cdot \lambda_2(\Ar)
\ll \sum_{\frac{q}{d} \mid r}\lambda_2\left(A_{d\bm{p}} \cap A_{\frac{r}{q/d}\bm{p}}\right) - \lambda_2(A_{d\bm{p}}) \cdot \lambda_2\left(A_{\frac{r}{q/d}\bm{p}}\right).
\end{split}
\]

Since $\psi(m) = O(m^{-\varepsilon})$ and $\gamma$ is a non-Liouville irrational, we
satisfy the assumptions of \eqref{assumptions} for some $C,c,\eta,\varepsilon$.
 Writing $M := \left\lfloor \frac{\eta+1}{\varepsilon}\right\rfloor, K = (2Cc)^{1/\varepsilon} +1$, an application of Lemma \ref{overlap_estim_new} thus gives

\[
\lambda_2\left(A_{d\bm{p}} \cap A_{\frac{r}{q/d}\bm{p}}\right) = 0
\]
if $r \geq d^M K$, and
\[
\lambda_2\left(A_{d\bm{p}} \cap A_{\frac{r}{q/d}\bm{p}}\right) \leq \lambda_2\left(A_{d\bm{p}}\right)\lambda_2\left(A_{\frac{r}{q/d}\bm{p}}\right) + 
4\gcd\left(d,\frac{r}{q/d}\right)\frac{\psi(q)}{d}
\]
otherwise. Thus replacing $r' = \frac{r}{q/d}$, we deduce

\[\begin{split}
\sum_{d\mid q}\sum_{\substack{\bq, \in \mathbb{Z}^2\\ \lvert \bq\rvert = q\\
\gcd(\bq) = d}}
\sum_{r < q}
\sum_{\substack{\br \in \mathbb{Z}^2\\ \lvert \br\rvert = r\\ \br \parallel \bq}} \lambda_2(\Aq \cap \Ar)-\lambda_2(\Aq)\cdot \lambda_2(\Ar)
\ll
\sum_{d\mid q}\sum_{\substack{\bq, \in \mathbb{Z}^2\\ \lvert \bq\rvert = q\\
\gcd(\bq) = d}}
\sum_{\substack{r' < \min\left\{d, \frac{d^{M+1} K}{q}\right\}}}
\frac{\psi(q)}{d}\gcd\left(d,r'\right).
\end{split}
\]

Since the number of $\bq \in \mathbb{Z}^2$ with $\lvert \bq \rvert = q$ and $\gcd(\bq) = d$ is bounded above by $4\frac{q}{d}$, we can further bound the right hand side above by

\[
\begin{split}
\psi(q)q
\sum_{d\mid q} \frac{1}{d^2}
\sum_{r' \leq \min\{d, \frac{d^{M+1}K}{q}\}}
(d,r')
 &\leq \psi(q)q
\sum_{\substack{d\mid q\\
d \geq (q/K)^{1/(M+1)}}
}\frac{1}{d^2}
\sum_{r' \leq d}
(d,r').
\end{split}
\]
Note that $\sum_{r' \leq d}
(d,r') \leq d \tau(d)$ where $\tau$ denotes the number of divisors. Since $\tau(d) \leq \tau(q) \ll_{\varepsilon} q^{\varepsilon}$ for any $\varepsilon > 0$, we can further bound the above expression by

\[
\psi(q)q
\sum_{\substack{d\mid q\\
d \geq (q/K)^{1/(M+1)}}
}\frac{\tau(d)}{d} \leq 
\psi(q)q
\tau(q)^2
\frac{1}{(q/K)^{1/(M+1)}} \ll \psi(q)q^{1 - 1/(M+2)} \leq q\psi(q).
\]
Thus the statement follows after summing over $q \leq Q$.
    \end{proof}
    
Now we have all ingredients to conclude the proof of Lemma \ref{sneaky_order}:
        Let
        \[X_1(\bm{\alpha}) := \sum_{\substack{\bm{u} \preceq \bq\\ \lvert \bm{u} \rvert = \lvert \bq \rvert} } \mathds{1}_{A_{\bq}(\bq,\psi(\lvert \bq \rvert))}(\bm{\alpha}), \quad
        X_2(\bm{\alpha}) := \sum_{\substack{\lvert \bm{u} \rvert < \lvert \bq \rvert
        < \lvert \bm{v} \rvert} } \mathds{1}_{A_{\bq}(\bq,\psi(\lvert \bq \rvert))}(\bm{\alpha}), \quad
    X_3(\bm{\alpha}) := \sum_{\substack{\bq \preceq \bm{v} \\ \lvert \bm{v} \rvert = \lvert \bq \rvert} } \mathds{1}_{A_{\bq}(\bq,\psi(\lvert \bq \rvert))}(\bm{\alpha}),
        \]
        which can be interpreted as random variables with $\bm{\alpha} \stackrel{d}\sim Unif(\T^2)$. In this way, we can write

        \[ 
        \begin{split}\sum_{\bm{u} \preceq \bq,\br \preceq \bm{v}}
\lambda_2(\Aq \cap \Ar) 
- \left(\sum_{\bm{u} \preceq \bq,\br \preceq \bm{v}}\lambda_2(\Aq)\right)^2 &= \Var[X_1 + X_2 + X_3] 
\\&\leq 3 \left(\Var[X_1] +\Var[X_2]+ \Var[X_3] \right).
\end{split}
        \]
        Starting with $\Var[X_1]$ and $\Var[X_3]$, we observe that all vectors $\br,\bq$ with $\lvert \bq \rvert = \lvert \br \rvert $ are either non-parallel or satisfy $\br = \pm \bq$.
        Thus Proposition \ref{nonparallel} combined with an analogue of \eqref{diag_contr} shows
        \[\Var[X_1] \ll \sum_{\substack{\bm{u} \preceq \bq\\ \lvert \bm{u} \rvert = \lvert \bq \rvert}}
        \psi(\lvert\bq\rvert), \quad \Var[X_3] \ll \sum_{\substack{\bq \preceq \bm{v} \\ \lvert \bm{v} \rvert = \lvert \bq \rvert} } \psi(\lvert\bq\rvert).
        \]
    For $\Var[X_2]$, we apply Lemma \ref{overlap_estim_new} to $\tilde{\psi}(q) := \psi(q)\mathds{1}_{[\lvert \bm{u}\rvert+1,\lvert \bm{v}\rvert-1]}(q)$ to obtain

    \[\Var[X_2] \ll \sum_{\substack{\lvert \bm{u} \rvert < q
        < \lvert \bm{v} \rvert} } 
        q\psi(q) \ll  \sum_{\substack{\lvert \bm{u} \rvert < \bq
        < \lvert \bm{v} \rvert}}\psi(\lvert \bq\rvert),
        \]
        so combining the estimates above, the statement follows.

\subsection*{Further remarks}
    Here, we give a sketch that if $m = 2$, the non-Liouville condition can be dropped from Theorem \ref{main_thm_new}, when we still keep the assumption $\psi(q) = O(1/q^{\varepsilon})$:
 For $m \geq 2$, Lemma \ref{overlap_estim_new} is replaced by the estimate (note that $\lambda_m(\Aq) = 2^m\psi(\lvert \bm{q}\rvert )^m$)

\begin{equation}\label{overlap_highdim}
\lambda_m(\Aq \cap \Ar) \leq \left(
4\psi(q)\psi(r) + 4\frac{\psi(q)}{q}(q,r)\right)^m.
\end{equation}

This can be obtained by defining a $m$-dimensional weight-function as a product of $m$ 1-dimensional weight-functions of the form \eqref{weight_fct}. On each coordinate, we apply the general estimates from Lemma \ref{overlap_estim_new} and note that with the notations from Lemma \ref{overlap_estim_new}, $\frac{(d,e)}{e} = \frac{(q,r)}{q}$, which proves \eqref{overlap_highdim}.\\

 For $m = 2, n = 2$, we use $\psi(r) \ll \frac{1}{r^{\varepsilon}}$ to further bound the right-hand side of \eqref{overlap_highdim} by 
 \[
 \psi(q)^2\psi(r)^2 + O\left(\frac{\psi(q)^2}{q}(q,r)^{1-\varepsilon} + \frac{\psi(q)^2}{q^2}(q,r)^2\right).
 \]
 
 For fixed $r < q$ there are $\ll (r,q)^{n-1}$ many parallel pairs of vector $\bq,\br \in \mathbb{Z}^n$ with norms
 $q$ and $r$ respectively. Thus for fixed $q$ we obtain
 \begin{align*}
\sum_{\substack{\bq \in \mathbb{Z}^2\\ \lvert \bq\rvert = q}}
\sum_{r < q}\sum_{\substack{\br \in \mathbb{Z}^2\\ \lvert \br\rvert = r\\ \br \parallel \bq}} \lambda_2(\Aq \cap \Ar) - \lambda_2(\Aq) \lambda_2(\Ar)
\ll \frac{\psi(q)^2}{q}\sum_{r < q}(q,r)^{2-\varepsilon} + \frac{\psi(q)^2}{q^2}\sum_{r < q}(q,r)^3
\ll q\psi(q)^2,
 \end{align*}
 where we used
 $\tau(q) \ll_{\varepsilon} q^{\varepsilon}$ and $\sum_{r < q}(q,r)^k \ll q^k$ for $k \geq 3$.
 If $n \geq 3$ and $m \geq 1$, then $\psi(r) \leq 1/2$ and $\sum_{r < q}(q,r)^k \ll q^k, k \geq 3$ alone is enough to show that

\[
\sum_{\substack{\bq \in \mathbb{Z}^2\\ \lvert \bq\rvert = q}}
\sum_{r < q}\sum_{\substack{\br \in \mathbb{Z}^2\\ \lvert \br\rvert = r\\ \br \parallel \bq}} \lambda_n(\Aq \cap \Ar) - \lambda_n(\Aq) \lambda_n(\Ar) \ll q^{n-1}\psi(q)^m.
\]
The rest of the proof follows along the lines of the case $n = 2, m = 1$.

        \subsection*{Acknowledgements} 
        
MH was supported by the EPSRC grant EP/X030784/1. The author would like to thank Christoph Aistleitner and Victor Beresnevich for many valuable discussions and comments on an earlier version of this manuscript.

\end{document}